\documentclass[12pt]{article}

\usepackage{amstext}
\usepackage{graphicx}
\usepackage{caption}
\usepackage{subcaption}
\usepackage{multirow}

\title{Dynamics of Two Coupled van der Pol Oscillators with Delay Coupling Revisited}
\author{Mark Gluzman\\
	Center for Applied Mathematics\\ Cornell University\\and\\
	 Richard Rand\\
	 Dept. Mathematics and Dept. Mechanical and Aerospace Engineering\\ Cornell University}

\begin{document}

\date{}
\maketitle

\vspace*{1 in}

\section{Abstract}
The problem of two van der Pol oscillators coupled by velocity delay terms was studied by Wirkus and Rand in 2002 \cite{wirkus}.
The small-$\varepsilon$ analysis resulted in a slow flow which contained delay terms.  To simplify the analysis, Wirkus and Rand followed
a common procedure of replacing the delay terms by non-delayed terms, a step said to be valid for small $\varepsilon$, resulting in a slow flow which
was an ODE rather than a DDE (delay-differential equation).  In the present paper we consider the same problem but leave the delay terms in the slow flow,
thereby offering an evaluation of the approximate simplification made in \cite{wirkus}.\\

\newpage

\section{Introduction}

In a recent paper by Sah and Rand \cite{simo}, it was shown that a class of nonlinear oscillators with delayed self-feedback gave rise via a perturbation solution to a DDE slow flow.  An exact solution was obtained to the DDE slow flow and was compared to the approximate solution resulting from the common approach of replacing the delayed variables in the slow flow with non-delayed variables, thereby replacing the DDE slow flow with an easier to solve ODE slow flow.\\

The paper by Sah and Rand \cite{simo} was motivated by numerous papers in the literature of nonlinear dynamics in which the DDE slow flow is replaced by an approximate ODE slow flow, for example  \cite{atay}, \cite{tina}, \cite{wirkus}.\\

In particular, in 2002 Wirkus and Rand wrote a paper entitled ``The Dynamics of Two Coupled van der Pol Oscillators with Delay Coupling" \cite{wirkus} in which averaging was used to derive a slow flow which governed the stability of the in-phase mode.  In studying the resulting slow flow, Wirkus and Rand replaced certain delayed quantities with non-delayed versions in order to simplify the analysis.  In this paper we reexamine the previously studied system with the idea of leaving the delayed quantities in the slow flow.\\

\section{Delay-Coupled van der Pol Oscillators}

The governing equations are:
\begin{eqnarray}
\label{goo1}
\ddot x_1+x_1-\varepsilon(1-x_1^2)\dot x_1=\varepsilon \alpha \dot x_2(t-T)\\
\label{goo2}
\ddot x_2+x_2-\varepsilon(1-x_2^2)\dot x_2=\varepsilon \alpha \dot x_1(t-T)
\end{eqnarray}
This system admits an in-phase mode in which $x_1=x_2=y(t)$ where
\begin{eqnarray}
\ddot y+y-\varepsilon(1-y^2)\dot y=\varepsilon \alpha \dot y(t-T)
\label{foo1}
\end{eqnarray}

In order to investigate the stability of the in-phase mode $y(t)$, we will first obtain an approximate expression for it using Lindstedt's method.  We replace independent variable $t$ by stretched time $\tau$:
\begin{eqnarray}
\tau = \omega t, \mbox{ ~~~~~~ where }\omega = 1 + \varepsilon k + O(\varepsilon^2)
\end{eqnarray}
so that eq.(\ref{foo1}) becomes:

\begin{eqnarray}
 \omega^2 y''+y-\varepsilon(1-y^2)\omega y'=\varepsilon \alpha  \omega y'(\tau-\omega T)
\label{foo2}
\end{eqnarray}
where primes represent differentiation with respect to $\tau$.  We expand $y$ in a power series in $\varepsilon$,
\begin{eqnarray}
y=y_0 + \varepsilon y_1 + \cdots
\label{foo3}
\end{eqnarray}
and substitute (\ref{foo3}) into (\ref{foo2}).  After collecting terms, we get
\begin{eqnarray}
\label{foo4}
 y_0''+ y_0 &=& 0\\
 \label{foo5}
  y_1''+ y_1 &=& -2 k y_0'' +(1-y_0^2)y_0' +\alpha y_0'(\tau-T)
\end{eqnarray}
We take the solution of (\ref{foo4}) to be
\begin{eqnarray}
y_0 = R \cos \tau
\label{foo6}
\end{eqnarray}
and substitute (\ref{foo6}) into (\ref{foo5}), giving
\begin{eqnarray}
\label{foo7}
y_1''+ y_1 &=& (2 k R+\alpha R \sin T) \cos\tau + (-R + \frac{R^3}{4}- \alpha R \cos T)\sin \tau
                +\frac{R^3}{4} \sin(3\tau)
\end{eqnarray}
Removing secular terms in (\ref{foo7}), we obtain
\begin{eqnarray}
\label{foo8}
k=-\frac{\alpha}{2}\sin T ~~~~~~~~~\mbox{and}~~~~~~ R=2\sqrt{1+\alpha\cos T}
\end{eqnarray}
Thus the in-phase mode $x_1=x_2=y(t)$, the stability of which is desired, is given by the approximation
\begin{eqnarray}
\label{foo9}
 y(t) = R \cos \omega t = 2\sqrt{1+\alpha\cos T} \cos \left(1-\frac{\alpha}{2}\varepsilon\sin T\right) t
\end{eqnarray}

\section{Stability of the In-Phase Mode}

In order to study the stability of the motion (\ref{foo9}), we set
\begin{eqnarray}
\label{foo10}
x_1 = y(t) + w_1 ~~~~~~~\mbox{and}~~~~~x_2 = y(t) + w_2
\end{eqnarray}
where $w_1$ and $w_2$ are deviations off of the in-phase mode $x_1=x_2=y(t)$.  Substituting (\ref{foo10}) into
(\ref{goo1}) and (\ref{goo2}) and using (\ref{foo1}), we obtain the following equations on  $w_1$ and $w_2$, linearized about
$w_1=w_2=0$:
\begin{eqnarray}
\label{goo3}
 \ddot w_1+(1+2\varepsilon y \dot y)w_1-\varepsilon(1-y^2) \dot w_1=\varepsilon \alpha  \dot w_2(t-T)\\
\label{goo4}
\ddot  w_2+(1+2\varepsilon y \dot y)w_2-\varepsilon(1-y^2) \dot w_2=\varepsilon \alpha  \dot w_1(t-T)
\end{eqnarray}
Equations (\ref{goo3}) and (\ref{goo4}) can be uncoupled by setting
\begin{eqnarray}
\label{goo5}
z_1 = w_1 + w_2 ~~~~~~~\mbox{and}~~~~~z_2 = w_1 - w_2
\end{eqnarray}
giving
\begin{eqnarray}
\label{goo6}
\ddot z_1+(1+2\varepsilon y \dot y)z_1-\varepsilon(1-y^2) \dot z_1=\varepsilon \alpha  \dot z_1(t-T)\\
\label{goo7}
\ddot z_2+(1+2\varepsilon y \dot y)z_2-\varepsilon(1-y^2) \dot z_2=-\varepsilon \alpha  \dot z_2(t-T)
\end{eqnarray}
The only difference between these two equations is the sign of the right-hand side, which may be absorbed into a new coefficient, call it $\beta$, which equals either 1 or -1:
\begin{eqnarray}
\label{goo8}
\ddot u+(1+2\varepsilon y \dot y)u-\varepsilon(1-y^2) \dot u=\varepsilon \beta \alpha  \dot u(t-T)
\end{eqnarray}
where $u=z_1$ for $\beta=1$ and $u=z_2$ for $\beta=-1$.  In eq.(\ref{goo8}), $y(t)$ is given by eq.(\ref{foo9}).
In order to study the boundedness of solutions to eq.(\ref{goo8}), we return to using $\tau = \omega t$ as independent variable,
where $\omega=1-\frac{\alpha}{2}\varepsilon\sin T$:
\begin{eqnarray}
\label{goo9}
\omega^2 u''+(1+2\varepsilon \omega y  y')u-\varepsilon(1-y^2) \omega u'=\varepsilon \beta \alpha  \omega u'(\tau-\omega T)
\end{eqnarray}
We study (\ref{goo9}) by using the two variable perturbation method \cite{rhr}.  We let $\xi=\tau$ and
$\eta = \varepsilon \tau$, giving:
\begin{eqnarray}
\label{goo10}
\omega^2 (u_{\xi\xi}+2 u_{\xi\eta})+(1+2\varepsilon \omega y  y_\xi )u-\varepsilon(1-y^2) \omega u_\xi=\varepsilon \beta \alpha  \omega u_\xi(\xi-\omega T,\eta-\varepsilon \omega T)
\end{eqnarray}
where $y(\xi)=R\cos\xi$, $R=2\sqrt{1+\alpha\cos T}$,  $\omega=1-\frac{\alpha}{2}\varepsilon\sin T$ and where we have neglected terms of  $O(\varepsilon^2)$.  Now we expand $u = u_0 + \varepsilon u_1 + O(\varepsilon^2)$ and collect terms, giving:
\begin{eqnarray}
\label{goo11}
u_{0_{\xi\xi}}+ u_0=0~~~~~~~~~~~~~~~~~~~~~~~~~~~~~~~~~~~~~~~~~~~~~~~~~~~~~~~~~~~~~~~~~~
\end{eqnarray}
\begin{eqnarray}
\label{goo12}
\nonumber
u_{1_{\xi\xi}}+ u_1=-2 u_{0{_\xi\eta}}	-\alpha\sin T~ u_{0_{\xi\xi}}+8(1+\alpha\cos T)\cos \xi \sin \xi~ u_0&&\\
-4(1+\alpha\cos T)\cos^2\xi ~u0_{\xi} +\alpha\beta u_{0_{\xi}}(\xi-T,\eta-\varepsilon T)
\end{eqnarray}
We take the solution of (\ref{goo11}) in the form
\begin{eqnarray}
\label{goo13}
u_0=A(\eta)\cos\xi+B(\eta)\sin\xi																									
\end{eqnarray}
Note that $u_{0}(\xi-T,\eta-\varepsilon T)=A_d \cos(\xi-T)+B_d\sin(\xi-T)$, where
\begin{eqnarray}
\label{goo14}
A_d=A(\eta-\varepsilon T)   \mbox{~~~and~~~}  B_d=B(\eta-\varepsilon T)
\end{eqnarray}
Substituting (\ref{goo13}),(\ref{goo14}) into (\ref{goo12}) and eliminating secular terms gives the slow flow
\begin{equation}
\label{goo15}
{{d\,A}\over{d\,\eta}}=-A
-{{3\,\alpha\,A\,\cos T}\over{2}}
+{{\alpha\,B\,\sin T}\over{2}}
+{{\alpha\,{\it A_d}\,\beta\,\cos T}\over{2}}
-{{\alpha\,\beta\,{\it B_d}\,\sin T}\over{2}}
\end{equation}
\begin{equation}
\label{goo16}
{{d\,B}\over{d\,\eta}}=
-{{\alpha\,A\,\sin T}\over{2}}
-{{\alpha\,B\,\cos T}\over{2}}
+{{\alpha\,{\it A_d}\,\beta\,\sin T}\over{2}}
+{{\alpha\,\beta\,{\it B_d}\,\cos T}\over{2}}
\end{equation}
The slow flow (\ref{goo15}),(\ref{goo16}) is a system of linear delay-differential equations (DDEs).
A common approach to treating such a system is to replace the delayed variables by non-delayed variables, as in $A_d=A$ and $B_d=B$.
To support such a step, it is often argued that since a Taylor expansion gives $A(\eta-\varepsilon T)= A(\eta) + O(\varepsilon)$, the replacement of $A_d$ by $A$ is an approximation valid for small $\varepsilon$ \cite{atay},\cite{tina},\cite{simo}.  Let us follow this procedure and see what we get, and then compare results with what we would get by treating the system as a DDE. \\

Replacing $A_d$ and $B_d$ by $A$ and $B$, we obtain the ODE system:
\begin{equation}
\label{goo17}
\frac{d}{d\eta}
\left(
\begin{array}{c}
A \\ B
\end{array}\right)=
\left( \begin{array}{cc}
-1+\frac{\alpha}{2}(\beta-3)\cos T & \frac{\alpha}{2}(1-\beta)\sin T  \\
-\frac{\alpha}{2}(1-\beta)\sin T & -\frac{\alpha}{2}(1-\beta)\cos T  \end{array} \right)
\left(
\begin{array}{c}
A \\ B
\end{array}\right)
\end{equation}
Recall from eqs.(\ref{goo6})-(\ref{goo8}) that there are two cases, $\beta=+1$ and $\beta=-1$.
In the case that $\beta=-1$, the system (\ref{goo17}) becomes:
\begin{equation}
\label{goo18}
\frac{d}{d\eta}
\left(
\begin{array}{c}
A \\ B
\end{array}\right)=
\left( \begin{array}{cc}
-1-2\alpha\cos T & \alpha\sin T  \\
-\alpha\sin T & -\alpha\cos T  \end{array} \right)
\left(
\begin{array}{c}
A \\ B
\end{array}\right)
\end{equation}
This system of ODEs exhibits both Hopf and saddle-node bifurcations.  The Hopf bifurcations occur when the trace = $-1-3\alpha\cos T$ = 0 with positive determinant, i.e. when
\begin{equation}
\label{hoo1}
\mbox{Hopf bifurcations:}~~~~~~\alpha = -\frac{1}{3 \cos T}
\end{equation}
The saddle-node bifurcations occur when the determinant = $\alpha^2+\alpha \cos T+\alpha^2 \cos^2 T$ = 0, i.e. when
\begin{equation}
\label{hoo2}
\mbox{saddle-node bifurcations:}~~~~~~\alpha = 0 \mbox{~~~and~~~} \alpha = -\frac{\cos T}{1+ \cos^2 T}
\end{equation}
Now let us consider the other case, $\beta = +1$.  In this case the ODE system (\ref{goo17}) becomes:
\begin{equation}
\label{hoo3}
\frac{d}{d\eta}
\left(
\begin{array}{c}
A \\ B
\end{array}\right)=
\left( \begin{array}{cc}
-1-\alpha\cos T & 0  \\
0 & 0  \end{array} \right)
\left(
\begin{array}{c}
A \\ B
\end{array}\right)
\end{equation}
From eq.(\ref{foo8}) we see that the amplitude of the in-phase mode, whose stability we are investigating, is given by
  $R=2\sqrt{1+\alpha\cos T}$.  Inspection of (\ref{hoo3}) shows that that system exhibits the birth of the in-phase mode when
\begin{equation}
\label{hoo4}
\mbox{birth of the in-phase mode:}~~~~~~\alpha = -\frac{1}{\cos T}
\end{equation}
We note that all of the bifurcations (\ref{hoo1}),(\ref{hoo2}),(\ref{hoo4}) were observed by Wirkus and Rand \cite{wirkus} in
their original work on this system.\\

In what follows, we return to the DDE system (\ref{goo15})-(\ref{goo16}), but we do not make the simplifying assumption of replacing $A_d$ by $A$ and $B_d$ by $B$.  In order to treat the DDE system, we set:
\begin{equation}
\label{hoo5}
A=P e^{\lambda \eta},~~~B=Q e^{\lambda \eta}, ~~~A_d=P  e^{\lambda (\eta-\varepsilon T)},~~~B_d=Q  e^{\lambda (\eta-\varepsilon T)}
\end{equation}
where $P$ and $Q$ are constants.  We are particularly interested in the effect of the DDE slow flow on Hopf bifurcations, eq.(\ref{hoo1}).  Thus we restrict ourselves to the case $\beta=-1$, whereby we obtain the following pair of algebraic equations on $P$ and $Q$:
\begin{eqnarray}
\label{zoo1}
\pmatrix{-{{\alpha\,e^ {- {\it \lambda}\,{\it \varepsilon}\,T }\,\cos T
		}\over{2}}-{{3\,\alpha\,\cos T}\over{2}}-{\it \lambda}-1&{{\alpha\,e
		^ {- {\it \lambda}\,{\it \varepsilon}\,T }\,\sin T}\over{2}}+{{\alpha\,
		\sin T}\over{2}}\cr -{{\alpha\,e^ {- {\it \lambda}\,{\it \varepsilon}\,T
		}\,\sin T}\over{2}}-{{\alpha\,\sin T}\over{2}}&-{{\alpha\,e^ {-
		{\it \lambda}\,{\it \varepsilon}\,T }\,\cos T}\over{2}}-{{\alpha\,\cos T
}\over{2}}-{\it \lambda}\cr }
\left(
\begin{array}{c}
P \\ Q
\end{array}\right)=\left(
\begin{array}{c}
0 \\ 0
\end{array}\right)
\end{eqnarray}
For a nontrivial solution, the determinant must vanish:
\begin{eqnarray}
\label{zoo2}
\nonumber
\alpha\,\cos T\,\lambda\,e^ {- {\it \varepsilon}\,T\,\lambda }-{{
		\alpha^2\,\sin ^2T\,e^ {- {\it \varepsilon}\,T\,\lambda }}\over{2}}+
{{\alpha\,\cos T\,e^ {- {\it \varepsilon}\,T\,\lambda }}\over{2}}+
\alpha^2\,e^ {- {\it \varepsilon}\,T\,\lambda }~~~~~~~\\
+{{\alpha^2\,e^ {- 2\,
			{\it \varepsilon}\,T\,\lambda }}\over{4}}+\lambda^2+2\,\alpha\,\cos T
\,\lambda+\lambda-{{\alpha^2\,\sin ^2T}\over{2}}+{{\alpha\,\cos T
	}\over{2}}+{{3\,\alpha^2}\over{4}}=0
\end{eqnarray}
For $\varepsilon$=0 and $\beta=-1$ the DDE system (\ref{goo15}),(\ref{goo16}) reduces to the ODE system (\ref{goo18}),
which exhibits a Hopf bifurcation (\ref{hoo1}) when $\alpha = -\frac{1}{3 \cos T}$.
By determining the location of the associated Hopf bifurcation in (\ref{zoo2}) for  $\varepsilon >0$ we may assess the accuracy of
the often made approximation based on replacing the delayed variables in the slow flow with non-delayed variables.\\

In the case of the ODE system (\ref{goo18}),  we obtained the conditions for a Hopf bifurcation, namely that there be a pair of pure imaginary eigenvalues, by requiring the trace of the matrix (\ref{goo18}) to be zero. In the case of the DDE system (\ref{goo15}),(\ref{goo16}), we seek a Hopf bifurcation by setting $\lambda=i \Omega$ in (\ref{zoo2}).\\

We seek a solution to (\ref{zoo2}) in the form of a perturbation series in $\varepsilon$ by expanding
\begin{eqnarray}
\label{zoo3}
T=T_0+\varepsilon T_1+\varepsilon^2 T_2 +\cdots\\
\label{zoo4}
\Omega=\Omega_0+\varepsilon \Omega_1+\varepsilon^2 \Omega_2 +\cdots
\end{eqnarray}
Separating (\ref{zoo2}) into real and imaginary parts and collecting like powers of $\varepsilon$ allows us to obtain the following expressions:
\begin{eqnarray}
\label{zoo5}
\cos T_0&=&-\frac{1}{3\alpha}\\
\label{zoo6}
\Omega_0&=&{{\sqrt{9\,\alpha^2-2}}\over{3}} \\
\label{zoo7}
T_1&=&-{{\sqrt{9\,\alpha^2-1}\,{\it T_0}}\over{9}}\\
\label{zoo8}
\Omega_1&=&-{{\left(18\,\alpha^2-5\right)\,{\it T_0}}\over{54\,\sqrt{9\,\alpha^2-2}}}\\
\label{zoo9}
T_2&=&{{\sqrt{9\,\alpha^2-1}\,\left(27\,\alpha^2-6\right)\,
		{\it T_0}^2+\left(162\,\alpha^4-36\,\alpha^2+2\right)\,{\it T_0}
	}\over{1458\,\alpha^2-162}}\\
\label{zoo10}
\nonumber
\Omega_2&=&{{\sqrt{9\,\alpha^2-2}\,\left(\left(-8019\,\alpha^6+5346
		\,\alpha^4-1206\,\alpha^2+91\right)\,{\it T_0}^2
				+\sqrt{9\,\alpha^2-1
		}\,\left(648\,\alpha^4-324\,\alpha^2+40\right)\,{\it T_0}\right)
	}\over{157464\,\alpha^4-69984\,\alpha^2+7776}}\\
\end{eqnarray}

\newpage

\section {Results and Conclusions}

To summarize our treatment of the original coupled van der Pol eqs.(\ref{goo1}),(\ref{goo2}), we first used Lindstedt's method to obtain an approximate expression for the in-phase mode, eq.(\ref{foo9}).  Then we studied the stability of the in-phase mode by applying the two variable perturbation method to eq.(\ref{goo9}).  This resulted in the DDE slow flow (\ref{goo15}),(\ref{goo16}).  We investigated this system of equations in two ways:\\

1) First we followed a number of other works \cite{atay},\cite{tina},\cite{simo} by replacing the delayed variables  $A_d$ and $B_d$ by non-delayed variables $A$ and $B$.  This resulted in a system of ODEs (\ref{goo17}) which possessed Hopf and saddle-node bifurcations, in agreement with the earlier work of Wirkus and Rand \cite{wirkus}.\\

2) Then we treated the slow flow system  (\ref{goo15}),(\ref{goo16}) as DDEs which resulted in a transcendental characteristic equation (\ref{zoo2}) which is harder to solve than the more familiar polynomial characteristic equations of ODEs.  We sought a series solution (\ref{zoo3}),(\ref{zoo4}) and obtained the results listed in eqs.(\ref{zoo5})-(\ref{zoo10}).\\

The results are plotted in Fig.1 where we show the critical delay for Hopf bifurcation $T_{Hopf}$ versus coupling strength $\alpha$ for $\varepsilon = 0.5$. The dashed curve is the analytical approximation based on replacing the delay terms in the slow flow (\ref{goo15}),(\ref{goo16}) by non-delayed variables, as in $A_d=A$ and $B_d=B$.  Thus the dashed curve corresponds to $\varepsilon=0$.  The dash-dot curve and the solid curve correspond to the analytical approximations respectively given by 2- and 3-term truncations of eq.(\ref{zoo3}).  The $+$ signs represent stability transitions obtained by numerical integration of the DDEs (\ref{goo15})-(\ref{goo16}).\\

The comparison between the various approximations for the critical delay for Hopf bifurcation shown in Fig.1 is further explored in Table 1, where we list the errors obtained using 1-, 2- and 3-term truncations of eq.(\ref{zoo3}) compared to values obtained by numerical integration of the DDEs (\ref{goo15})-(\ref{goo16}).  The maximum error is computed over the set $\alpha\in[\frac{\sqrt{2}}{3},~1]$. (Here $\alpha =\frac{\sqrt{2}}{3} $ is chosen because $(\alpha, T) = (\frac{\sqrt{2}}{3}, \frac{3\pi}{4})$ is a point of intersection of the bifurcation curves $\alpha = - \frac{cosT}{1+\cos^2 T}$ and $\alpha = - \frac{1}{3\cos T}$, compare eqs.(\ref{hoo1}),(\ref{hoo2}).) 
We consider \emph{absolute error} $\max\limits_{\alpha\in[\frac{\sqrt{2}}{3},~1]} |T(\alpha) - T_n(\alpha)|$, \emph{relative error} $\max\limits_{\alpha\in[\frac{\sqrt{2}}{3},~1]} |\frac{T(\alpha) - T_n(\alpha)}{T(\alpha)}|$, and \emph{percent error} $\max\limits_{\alpha\in[\frac{\sqrt{2}}{3},~1]} |\frac{T(\alpha) - T_n(\alpha)}{T(\alpha)}| \times 100\%$ , where $T(\alpha)$ is the value we get by numerical integration of the DDEs (\ref{goo15})-(\ref{goo16}) and $T_n(\alpha)$ is the value given by the n-term truncation of eq.(\ref{zoo3}). We also note that the maximum in all three cases was attained at $\alpha=1$. As previously noted, replacing the delay terms in the slow flow (\ref{goo15}),(\ref{goo16}) by non-delayed variables means that we take only the first term in eq.(\ref{zoo3}) and $n=1$.
In this case, the error incurred by omitting the delay terms in the slow flow is generally about 3 to 15$\%$ for small values of $\varepsilon$. On the other hand, the 3-term truncation of eq.(\ref{zoo3}) typically has the percent error less than 1\%.\\

\newpage

\begin{figure}[h!]
	\includegraphics[width=1\textwidth]{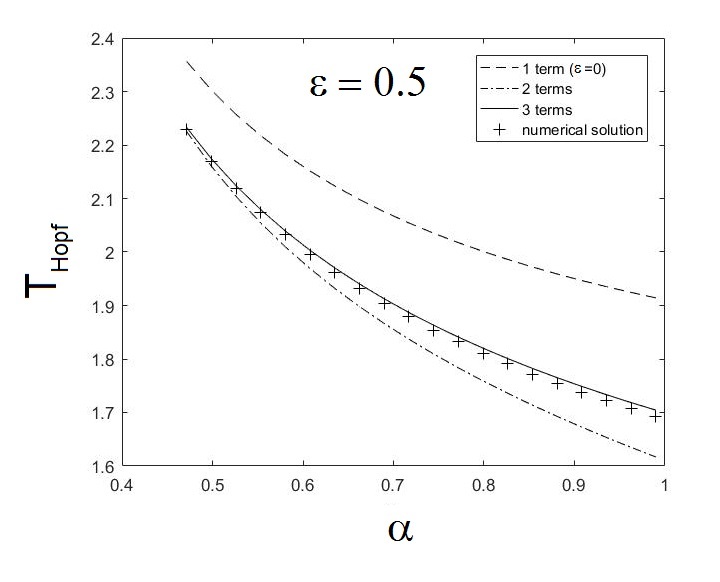}
	\caption{Critical delay for Hopf bifurcation versus $\alpha$ for $\varepsilon = 0.5$. Dashed curve is the analytical approximation based on replacing the delay terms in the slow flow (\ref{goo15}),(\ref{goo16}) by non-delayed variables, as in $A_d=A$ and $B_d=B$.  Thus the dashed curve corresponds to $\varepsilon=0$.  The dash-dot curve and the solid curve correspond to the analytical approximations respectively given by 2- and 3-term truncations of eq.(\ref{zoo3}).  The $+$ signs represent stability transitions obtained by numerical integration of the DDEs (\ref{goo15})-(\ref{goo16}).}
\end{figure}

\begin{table}[!htbp]
\centering

\label{table}
\begin{tabular}{l|l|l|l|l|l|l|l|l|l|}
\cline{2-10}
                                                              & \multicolumn{3}{c|}{$\varepsilon=0.1$}                                         & \multicolumn{3}{c|}{$\varepsilon=0.3$}                                         & \multicolumn{3}{c|}{$\varepsilon=0.5$}                                         \\ \hline
\multicolumn{1}{|l|}{n terms in (\ref{zoo3}) } & n=1                     & n=2                     & n=3                     & n=1                     & n=2                     & n=3                     & n=1                     & n=2                     & n=3                     \\ \hline
\multicolumn{1}{|l|}{\multirow{2}{*}{absolute error}}         & \multirow{2}{*}{0.056}  & \multirow{2}{*}{0.0025} & \multirow{2}{*}{0.001}  & \multirow{2}{*}{0.1499} & \multirow{2}{*}{0.0285} & \multirow{2}{*}{0.003}  & \multirow{2}{*}{0.2218} & \multirow{2}{*}{0.0756} & \multirow{2}{*}{0.012}  \\
\multicolumn{1}{|l|}{}                                        &                         &                         &                         &                         &                         &                         &                         &                         &                         \\ \hline
\multicolumn{1}{|l|}{\multirow{2}{*}{relative error}}         & \multirow{2}{*}{0.0307} & \multirow{2}{*}{0.0013} & \multirow{2}{*}{0.0005} & \multirow{2}{*}{0.085}  & \multirow{2}{*}{0.0162} & \multirow{2}{*}{0.0017} & \multirow{2}{*}{0.1311} & \multirow{2}{*}{0.0447} & \multirow{2}{*}{0.0071} \\
\multicolumn{1}{|l|}{}                                        &                         &                         &                         &                         &                         &                         &                         &                         &                         \\ \hline
\multicolumn{1}{|l|}{\multirow{2}{*}{percent error}}         & \multirow{2}{*}{3.07\%} & \multirow{2}{*}{0.13\%} & \multirow{2}{*}{0.05\%} & \multirow{2}{*}{8.5\%}  & \multirow{2}{*}{1.62\%} & \multirow{2}{*}{0.17\%} & \multirow{2}{*}{13.11\%} & \multirow{2}{*}{4.47\%} & \multirow{2}{*}{0.71\%} \\
\multicolumn{1}{|l|}{}                                        &                         &                         &                         &                         &                         &                         &                         &                         &                         \\ \hline
\end{tabular}
\caption{Errors in critical delay for Hopf bifurcation produced by 1-, 2- and 3-term truncations of eq.(\ref{zoo3}) as compared to values obtained by numerical integration of the slow flow  (\ref{goo15})-(\ref{goo16}).  The errors are defined as
	\emph{absolute error} $\max\limits_{\alpha\in[\frac{\sqrt{2}}{3},~1]} |T(\alpha) - T_n(\alpha)|$, \emph{relative error} $\max\limits_{\alpha\in[\frac{\sqrt{2}}{3},~1]} |\frac{T(\alpha) - T_n(\alpha)}{T(\alpha)}|$, and \emph{percent error} $\max\limits_{\alpha\in[\frac{\sqrt{2}}{3},~1]} |\frac{T(\alpha) - T_n(\alpha)}{T(\alpha)}| \times 100\%$ , where $T(\alpha)$ is the value we get by numerical integration of the DDEs (\ref{goo15})-(\ref{goo16}) and $T_n(\alpha)$ is the value given by the n-term truncation of eq.(\ref{zoo3}).
	}
\end{table}
\newpage
\vspace*{.1 in}

\section {Acknowledgement}

The authors thank Alex Bernstein for reading the manuscript.\\

\end{document}